\documentclass[12pt]{article}

\usepackage{amsmath}
\usepackage{amsfonts}
\usepackage{bigints}

\begin{document}

\renewcommand{\P}{\par \rm Proof:\ }
\newcommand{\Pe}{$\hfill{\Box}$\bigskip}
\newcommand{\Pes}{$\hfill{\Box}$}
\newcommand{\A}{\mbox{${{{\cal A}}}$}}

%1111111111111111111111111111111111111111111111111111111111111111111111111111

\author{Attila Losonczi}
\title{Means of unbounded sets}

\date{25 February 2018}

\newtheorem{thm}{\qquad Theorem}[section]
\newtheorem{prp}[thm]{\qquad Proposition}
\newtheorem{lem}[thm]{\qquad Lemma}
\newtheorem{cor}[thm]{\qquad Corollary}
\newtheorem{rem}[thm]{\qquad Remark}
\newtheorem{ex}[thm]{\qquad Example}
\newtheorem{df}[thm]{\qquad Definition}
\newtheorem{prb}{\qquad Problem}

\maketitle

\begin{abstract}

\noindent

We  study generalized means whose domain may contain unbounded sets as well. We investigate usual properties of this type of  means and also new attributes that regard for such means only.
We examine how a mean defined on bounded sets can be extended to this type of mean. We generalize some classic means and also present many new examples for means defined on unbounded sets.

\noindent
\footnotetext{\noindent
AMS (2010) Subject Classifications: 26E60, 40A05, 28A10 \\

} 

\end{abstract}

%-------------------------------------------------------------------------------------------------------------------------------------------Introduction----------------------------------
\section{Introduction}
This paper can be considered as a natural continuation of the investigations started in \cite{lamis} and \cite{lamisii} where we started to build the theory of means on infinite sets. An ordinary mean is for calculating the mean of two (or finitely many) numbers. This can be extended in many ways in order to get a more general concept where we have a mean on some infinite bounded subsets of $\mathbb{R}$. The various general properties of such means, the relations among those means were studied thoroughly in \cite{lamis} and \cite{lamisii}. 

In this paper our main aim is to study means that domain may contain unbounded sets as well.

First we investigate the general properties of such means. We also study already known properties for this type of means and new attributes are presented as well.  

Then we examine how a mean defined on bounded sets can be extended to a mean that is defined also on some unbounded sets. We check which properties of the original mean are inherited to the extension. 

We also present many new examples for means defined on unbounded sets too and we find natural generalizations for some classic means in order to get a mean defined on some unbounded sets as well. 

Finally we analyse the behavior of one of the most important generic means ${\cal{M}}^{\mu}\ (Avg)$ regarding unbounded measurable sets.

%------------------------------------------------------------------------------------------------------------------------------Basic notions and notations---------------------------
\subsection{Basic notions and notations}

Let us recall some very basic notions from \cite{lamis} and \cite{lamisii}.

\medskip

We call ${\cal{K}}$ an \textbf{ordinary mean} if it is for calculating the mean ${\cal{K}}(a_1,\dots,a_n)$ of finitely many numbers $a_1,\dots,a_n\in\mathbb{R}$.

\smallskip

A \textbf{generalized mean} is a function ${\cal{K}}:C\to \mathbb{R}$ where $C\subset P(\mathbb{R})$ consists of some (finite or infinite) subsets of $\mathbb{R}$ and $\inf H\leq {\cal{K}}(H)\leq\sup H$ holds for all $H\in C$.

\smallskip

A mean $\cal{K}$ is called \textbf{monotone} if $\sup H_1\leq\inf H_2$ implies that ${\cal{K}}(H_1)\leq {\cal{K}}(H_1\cup H_2)\leq {\cal{K}}(H_2)$. 
${\cal{K}}$ is \textbf{base-monotone} if $H_1,H_2\in Dom({\cal{K}}), H_1\cap H_2=\emptyset$ then $\min\{{\cal{K}}(H_1),{\cal{K}}(H_2)\}\leq{\cal{K}}(H_1\cup H_2)\leq\max\{{\cal{K}}(H_1),{\cal{K}}(H_2)\}.$

${\cal{K}}$ is \textbf{part-slice-continuous} if $H_1,H_2\in Dom({\cal{K}})$ then $H_2^{\varliminf H+},H_2^{\varlimsup H-}\in Dom({\cal{K}})$ and $f(x)={\cal{K}}(H_1\cup H_2^{x-})$ and $g(x)={\cal{K}}(H_1\cup H_2^{x+})$ are continuous where $Dom(f)=\{x:H_1\cup H_2^{x-}\in Dom({\cal{K}})\}, Dom(g)=\{x:H_1\cup H_2^{x+}\in Dom({\cal{K}})\}$.

$\cal{K}$ is \textbf{finite-independent} if $H$ being infinite implies that ${\cal{K}}(H)={\cal{K}}(H\cup V)={\cal{K}}(H-V)$ where $V$ is any finite set.

%${\cal{K}}$ is \textbf{part-shift-monotone} if $x>0,H_1\cap H_2=H_1\cap H_2+x=\emptyset$ then ${\cal{K}}(H_1\cup H_2)\leq {\cal{K}}(H_1\cup H_2+x)$.

%\begin{df}Let ${\cal{K}}$ be a mean, $H\in Dom({\cal{K}})$. We say that a bounded $V$ is small for $H$ if $H\cup (V+x)\in Dom({\cal{K}}),H-(V+x)\in Dom({\cal{K}}),{\cal{K}}(H\cup V+x)={\cal{K}}(H-(V+x))={\cal{K}}(H)\ \forall x\in\mathbb{R}\}$.  
%\end{df}

\medskip

Throughout this paper $\lambda$ will denote the Lebesgue measure.

\begin{df}(cf. \cite{lambm} Def 2.1) Let $\mu$ be a Borel measure on $\mathbb{R}$. Let $H\subset\mathbb{R}$ be a $\mu$-measurable set such that $0<\mu(H)<+\infty$. Then 
$${\cal{M}}^{\mu}(H)=\frac{\int\limits_Hx d\mu}{\mu(H)}$$
is a mean defined on unbounded subsets as well.
\end{df}

We get a special case for Hausdorff measures.

\begin{df}\label{davg}Let $\mu^s$ denote the s-dimensional Hausdorff measure ($0\leq s\leq 1$). If $0<\mu^s(H)<+\infty$ (i.e. $H$ is an $s$-set) and $H$ is $\mu^s$ measurable then $$Avg(H)=\frac{\int\limits_H x\ d\mu^s}{\mu^s(H)}.$$
If $0\leq s\leq 1$ then set $Avg^s=Avg|_{\{\text{measurable s-sets}\}}$. E.g. $Avg^1$ is $Avg$ on all Lebesgue measurable sets with positive finite measure.
\end{df}

If $H\subset\mathbb{R},x\in\mathbb{R}$ then set $H+x=\{h+x:h\in H\}$. Similarly $\alpha H=\{\alpha h:h\in H\}\ (\alpha\in\mathbb{R})$.
We use the convention that this operation $+$ has to be applied prior to the set theoretical operations, e.g. $H\cup K\cup L+x=H\cup K\cup (L+x)$.

\smallskip

The extended real line is $\bar{\mathbb{R}}=\mathbb{R}\cup\{-\infty,+\infty\}$ equipped with the usual topology: the neighbourhood of $+\infty$ is $\{(k,+\infty]:k\in\mathbb{R}\}$ and similarly $\{[-\infty,k):k\in\mathbb{R}\}$ for $-\infty$. 

For $K\subset\mathbb{R},\ y\in\bar{\mathbb{R}}$ let us use the notations $$K^{y-}=K\cap(-\infty,y],K^{y+}=K\cap[y,+\infty),K^{+\infty-}=K^{-\infty+}=K.$$

If $x>y$ then let $[x,y]$ denote the interval $[y,x]$.

\smallskip

Let us define some usual operations, relation with $\pm\infty$: $(+\infty)+(+\infty)=+\infty,\ (-\infty)+(-\infty)=-\infty$, if $r\in\mathbb{R}$ then $r+(+\infty)=+\infty,\ r+(-\infty)=-\infty,\ -\infty<r<+\infty$.

\smallskip

$cl(H), H'$ will denote the closure and accumulation points of $H\subset\mathbb{R}$ respectively. Let $\varliminf H=\inf H',\ \varlimsup H=\sup H'$ for $H\subset\mathbb{R},\ H'\ne\emptyset$.

\smallskip

Usually ${\cal{K}},{\cal{M}}$ will denote means, $Dom({\cal{K}})$ denotes the domain of ${\cal{K}}$.

%------------------------------------------------------------------------------------------------------------------------------Properties of means defined on unbounded sets---------------------------
\section{Properties of means defined on unbounded sets}

In the sequel a mean ${\cal{K}}$ is always a mean defined on some unbounded sets as well rather than bounded sets only.

\medskip

If a mean  is defined on some unbounded sets then we require the usual basic property, internality that is $$\inf H\leq {\cal{K}}(H)\leq\sup H.$$
However almost always we require the stroger condition, strong internality $\varliminf H\leq {\cal{K}}(H)\leq\varlimsup H$ when $H'\ne\emptyset$. 

\medskip

The properties of $Dom\ {\cal{K}}$ that we require are: 
$Dom\ {\cal{K}}$ must be closed under finite union, intersection and
if $H\in Dom\ {\cal{K}}, I$ is an interval (finite or infinite) then $H\cap I\in Dom({\cal{K}})$ if $H\cap I\ne\emptyset$.

\medskip

We got used to the fact regarding means on bounded sets that it can happen that $\forall h\in H\ h<{\cal{K}}(H)$ or the opposite way around $\forall h\in H\ h>{\cal{K}}(H)$ (when ${\cal{K}}(H)=\sup H\not\in H$ or ${\cal{K}}(H)=\inf H\not\in H$ respectively). The same scenario can occur on means on unbounded sets too i.e. ${\cal{K}}(H)$ can be either $+$ or $-\infty$ whenever $H\subset\mathbb{R}$ (i.e. $\pm\infty\not\in H$).

\begin{df}Let $H\in Dom\ {\cal{K}}$ be unbounded. We call $H$ essentially unbounded above regarding ${\cal{K}}$ if $\forall x\in\mathbb{R}\ \exists y\in\mathbb{R}$ such that ${\cal{K}}(H^{x-})<{\cal{K}}(H^{y-})$. We call $H$ essentially unbounded below regarding ${\cal{K}}$ if $\forall x\in\mathbb{R}\ \exists y\in\mathbb{R}$ such that ${\cal{K}}(H^{x+})>{\cal{K}}(H^{y+})$. We call $H$ essentially unbounded regarding ${\cal{K}}$ if it is essentially unbounded above and below.
\end{df}

Now we enumerate some properties of means defined on unbounded sets as well that we refer and analyze to later. 

\begin{itemize}

\item ${\cal{K}}$ is \textbf{i-strong-internal} if $H\in Dom({\cal{K}}),\ H'\ne\emptyset$ then $\inf (H'-\{-\infty\})\leq {\cal{K}}(H)\leq\sup (H'-\{+\infty\})$.

\item ${\cal{K}}$ is \textbf{slice-continuous} if $H\in Dom({\cal{K}})$ then $\forall y\in\mathbb{R}\ f_y(x)={\cal{K}}(H\cap[x,y])$ is continuous on the extended real line $\bar{\mathbb{R}}$ where $Dom( f_y)=\{x:H\cap[x,y]\ne\emptyset\}$ (cf. \cite{lamisii} 2.5).

\item ${\cal{K}}$ is \textbf{i-slice-continuous} if $H\in Dom({\cal{K}})$ then $f(x)={\cal{K}}(H^{x-})$ and $g(x)={\cal{K}}(H^{x+})$ are continuous on the extended real line $\bar{\mathbb{R}}$ where $Dom( f)=\{x:H^{x-}\ne\emptyset\},\ Dom(g)=\{x:H^{x+}\ne\emptyset\}$.

\item \textbf{The bounded sets are small for sets with infinite mean} if $K\in Dom\ {\cal{K}}$  is bounded and ${\cal{K}}(H)=+\infty$ then ${\cal{K}}(H\cup K)=+\infty$ and similarly 
if ${\cal{K}}(H)=-\infty$ then ${\cal{K}}(H\cup K)=-\infty$.

\item $K\in Dom\ {\cal{K}}$ is said to be t-infinite regarding ${\cal{K}}$ if $H\subset\mathbb{R}$ is bounded and $H\cup K+x\in Dom\ {\cal{K}}$ then $\lim\limits_{x\to+\infty}{\cal{K}}(H\cup K+x)=+\infty,\ \lim\limits_{x\to-\infty}{\cal{K}}(H\cup K+x)=-\infty$.
${\cal{K}}$ is \textbf{interval-infinite} if $I$ is a non-degenerative finite interval then it is t-infinite regarding ${\cal{K}}$.

\item $K\in Dom\ {\cal{K}}$ is said to be t-continuous regarding ${\cal{K}}$ if $H,H\cup K+x\in Dom\ {\cal{K}}$ then the function $x\mapsto{\cal{K}}(H\cup K+x)$ is continuous.
${\cal{K}}$ is called \textbf{interval-continuous} if $I$ is a non-degenerative finite interval then it is t-continuous regarding ${\cal{K}}$.

\item ${\cal{K}}$ is called \textbf{finite} if ${\cal{K}}(H)$ is finite for all $H\in Dom\ {\cal{K}}$.

\item ${\cal{K}}$ is called \textbf{subset-finite} if $H,K\in Dom\ {\cal{K}},\ |{\cal{K}}(H)|<+\infty,\ K\subset H$ then $|{\cal{K}}(K)|<+\infty$.

\item ${\cal{K}}$ is called \textbf{bounded-finite} if $|{\cal{K}}(H)|<+\infty,\ K\in Dom\ {\cal{K}}$ is bounded then $|{\cal{K}}(H\cup K)|<+\infty$.

\item $H\in Dom\ {\cal{K}}$ is called limit-finite for ${\cal{K}}$ if $\lim\limits_{x\to+\infty}{\cal{K}}(H^{x+})-\inf H^{x+}=\lim\limits_{x\to-\infty}{\cal{K}}(H^{x-})-\sup H^{x-}=0$. ${\cal{K}}$ is called \textbf{limit-finite} if all $H\in Dom\ {\cal{K}},\ |{\cal{K}}(H)|<+\infty$ are limit-finite. 

\item ${\cal{K}}$ is called \textbf{strong-base-monotone} if it is base-monotone and the following holds. Let $H_1,H_2,K\in Dom\ {\cal{K}}$ be bounded sets such that $H_1\subset H_2,\ H_2\cap K=\emptyset$. By base-monotonicity there are constants $0\leq c,d,c',d'\leq 1$ such that $c+d=c'+d'=1$ and 

${\cal{K}}(H_1\cup K)=c\cdot {\cal{K}}(H_1)+d\cdot {\cal{K}}(K),\ {\cal{K}}(H_2\cup K)=c'\cdot {\cal{K}}(H_2)+d'\cdot {\cal{K}}(K).$

If ${\cal{K}}(H_1\cup K)={\cal{K}}(H_1)={\cal{K}}(K)$ then let $c=0$.

Strong-base-monotonicity requires that $c\leq c'$ has to hold.

\end{itemize}

\begin{prp}Let ${\cal{K}}$ be i-slice-continuous and let the bounded sets are small for sets with infinite mean. Let $H$ be unbounded such that $H^{0-},H^{0+}\in Dom\ {\cal{K}}$ and ${\cal{K}}(H^{0-})=-\infty,\  {\cal{K}}(H^{0+})=+\infty$. Then $H\not\in Dom\ {\cal{K}}$.
\end{prp}
\P If assuming the contrary $H$ was in $Dom\ {\cal{K}}$ then ${\cal{K}}(H)={\cal{K}}(H^{-\infty+})=\lim\limits_{x\to-\infty}{\cal{K}}(H^{x+})$ by i-slice-continuity. But $x<0$ implies that ${\cal{K}}(H^{x+})={\cal{K}}(H^{0+})$ since $H^{x+}=[x,0)\cup H^{0+}$ and $[x,0)$ does not change the mean. Therefore ${\cal{K}}(H)={\cal{K}}(H^{0+})=+\infty$. Exactly the same way we would get that ${\cal{K}}(H)={\cal{K}}(H^{0-})=-\infty$ that is a contradiction.
\Pes

\begin{prp}\label{pafaiui}Let ${\cal{K}}$ be i-slice-continuous and let the bounded sets are small for sets with infinite mean. Let $H\in Dom\ {\cal{K}}$ be unbounded such that ${\cal{K}}(H^{a-})>-\infty,\  {\cal{K}}(H^{a+})=+\infty\ (a\in\mathbb{R})$. Then ${\cal{K}}(H)=+\infty$.
\end{prp}
\P ${\cal{K}}(H)={\cal{K}}(H^{-\infty+})=\lim\limits_{x\to-\infty}{\cal{K}}(H^{x+})$ by i-slice-continuity. But $x<a$ implies that ${\cal{K}}(H^{x+})={\cal{K}}(H^{a+})$ since $H^{x+}=[x,a)\cup H^{a+}$ and $[x,a)$ does not change the mean. Therefore ${\cal{K}}(H)={\cal{K}}(H^{a+})=+\infty$.
\Pes

\begin{prp}Let ${\cal{K}}$ be i-slice-continuous and let the bounded sets are small for sets with infinite mean. If $|{\cal{K}}(H)|<+\infty$ then $\forall x\in\mathbb{R}\  {\cal{K}}(H^{x+})<+\infty,\  {\cal{K}}(H^{x-})>-\infty$.
\end{prp}
\P Let us show the first and let $x\in\mathbb{R}$. Suppose indirectly that ${\cal{K}}(H^{x+})=+\infty$. Then by similar argument that we followed in the previous propositions one gets that ${\cal{K}}(H)=+\infty$ which is a contradiction.
\Pes

\begin{prp}Let ${\cal{K}}$ be slice-continuous, monotone and let the bounded sets are small for sets with infinite mean. Let $H\in Dom\ {\cal{K}}$ be unbounded such that ${\cal{K}}(H^{0-})=-\infty,\  {\cal{K}}(H^{0+})=+\infty$. Then for every $d\in\mathbb{R}$ there are sequences $(x_n),\ (y_n)$ such that $x_n\to-\infty,\ y_n\to+\infty$ and ${\cal{K}}(H\cap[x_n,y_n])=d$.
\end{prp}
\P Suppose $d>0$ (the remaining cases can be handled similarly).

Take $x_1<0$ such that $H\cap[x_1,0]\in Dom\ {\cal{K}}$. Then by slice-continuity one can find $y_1>0$ such that ${\cal{K}}(H\cap[x_1,y_1])=d$. If $x_{n-1},y_{n-1}$ is already given then take $x_n<x_{n-1}-1$ such that ${\cal{K}}(H\cap[x_n,y_{n-1}])<d$. Then by slice-continuity find $y_n>y_{n-1}$ such that ${\cal{K}}(H\cap[x_n,y_n])=d$. 

Clearly $(x_n)\to-\infty$, $(y_n)$ is increasing hence has a limit $y_n\to\beta\in\bar{\mathbb{R}}$. If $\beta<+\infty$ then we get ${\cal{K}}(H\cap[x_n,\beta])\to{\cal{K}}(H^{\beta-})$ but ${\cal{K}}(H^{\beta-})={\cal{K}}(H^{0-})=-\infty$. By monotonicity we have ${\cal{K}}(H\cap[x_n,y_n])<{\cal{K}}(H\cap[x_n,\beta])$ hence ${\cal{K}}(H\cap[x_n,y_n])\to-\infty$ - a contradiction.
\Pes

\begin{prp}Let ${\cal{K}}$ be i-slice-continuous. Then for every unbounded sets $H$ with finite mean there is a bounded set that is not small for $H$. 
\end{prp}
\P Assume the contrary and let $H$ be unbounded and ${\cal{K}}(H)=h\in\mathbb{R}$ such that all bounded sets are small for $H$. Then either $H^{h-}$ or $H^{h+}$ is unbounded. Say $H^{h+}$ is unbounded. Then ${\cal{K}}(H^{(h+1)+})>h$. Then ${\cal{K}}(H)={\cal{K}}(H^{-\infty+})=\lim\limits_{x\to-\infty}{\cal{K}}(H^{x+})$ by i-slice-continuity. But $x<h+1$ implies that ${\cal{K}}(H^{x+})={\cal{K}}(H^{(h+1)+})$ since $H^{x+}=[x,h+1)\cup H^{(h+1)+}$ and $[x,h+1)$ does not change the mean. Therefore $h={\cal{K}}(H)={\cal{K}}(H^{(h+1)+})>h$ that is a contradiction.
\Pes

%\begin{prp}Let the bounded sets be small for sets with infinite mean for a mean ${\cal{K}}$. If ${\cal{K}}(H)$ is infinite then ${\cal{K}}$ cannot be bi-slice-continuous (that is $f(x,y)={\cal{K}}(H^{y-}\cap H^{x+})$ is continuous on $\bar{\mathbb{R}}\times\bar{\mathbb{R}}$).
%\end{prp}
%\P Assuming the contrary one gets that ${\cal{K}}(H)=\lim\limits_{x\to+\infty}{\cal{K}}(H^{x-}\cap H^{x+})$ but all ${\cal{K}}(H^{x-}\cap H^{x+})$ equals to a same finite value e.g. ${\cal{K}}(H^{n-}\cap H^{n+})$ for some $n\in\mathbb{N}$ for which $H^{n-}\cap H^{n+}\ne\emptyset$.
%\Pes  ???????????

\begin{prp}\label{piscifhi}If ${\cal{K}}$ is i-slice-continuous, interval-infinite and the finite intervals are in $Dom\ {\cal{K}}$ then $\forall \epsilon>0$ there is $H\subset\mathbb{R}$ such that $\lambda(H)<\epsilon$ and ${\cal{K}}(H)=+\infty$.
\end{prp}
\P Let $I_1=[0,\frac{\epsilon}{3}]$. Then choose an interval $I_2$ such that $\lambda(I_2)<\frac{\epsilon}{2^2}$ and ${\cal{K}}(I_1\cup I_2)>2$. If we have choosen $I_1,\dots,I_{n-1}$ already then choose an interval $I_n$ such that $\lambda(I_n)<\frac{\epsilon}{2^n}$ and ${\cal{K}}(I_1\cup\dots\cup I_n)>n$.

Let $H=\bigcup\limits_1^{\infty}I_i$. Then obviously $\lambda(H)<\epsilon$ and by i-slice-continuity we get that ${\cal{K}}(H)=+\infty$.
\Pes

\begin{prp}If ${\cal{K}}$ is interval-infinite and the finite intervals are in $Dom\ {\cal{K}}$ then $\forall \epsilon>0$ there is $H\subset\mathbb{R}$ and sequences $(x_n),\ (y_n)$ such that $\lambda(H)<\epsilon,\ x_n\to-\infty,\ y_n\to+\infty$ and ${\cal{K}}(H\cap[x_n,y_n])$ is divergent.
\end{prp}
\P Let $I_1=[0,\frac{\epsilon}{3}]$. Then choose an interval $I_2=[a_{2},b_{2}]$ such that $\lambda(I_2)<\frac{\epsilon}{2^2}$ and ${\cal{K}}(I_1\cup I_2)>1$. Let $x_2=0$ and $y_2=b_2$.  Then choose an interval $I_3=[a_{3},b_{3}]$ such that $\lambda(I_3)<\frac{\epsilon}{2^3}$ and ${\cal{K}}(I_1\cup I_2\cup I_3)<-1$. Let $y_3=y_2$ and $x_3=a_3$.

If we have choosen $I_1,\dots,I_{2n-1}$ already then choose an interval $I_{2n}=[a_{2n},b_{2n}]$ such that $a_{2n}>b_{2n-2}+1,\ \lambda(I_{2n})<\frac{\epsilon}{2^{2n}}$ and ${\cal{K}}(I_1\cup\dots\cup I_{2n})>1$. Let $x_{2n}=x_{2n-1}$ and $y_{2n}=b_{2n}$.  Then choose an interval $I_{2n+1}=[a_{2n+1},b_{2n+1}]$ such that $b_{2n+1}<a_{2n-1}-1,\ \lambda(I_{2n+1})<\frac{\epsilon}{2^{2n+1}}$ and ${\cal{K}}(I_1\cup\dots\cup I_{2n+1})<-1$. Let $y_{2n+1}=y_{2n}$ and $x_{2n+1}=a_{2n+1}$.

Let $H=\bigcup\limits_1^{\infty}I_i$. Then obviously $x_n\to-\infty,\ y_n\to+\infty,\ \lambda(H)<\epsilon$ and ${\cal{K}}(H\cap[x_n,y_n])$ is divergent.
\Pes

\begin{prp}Let ${\cal{K}}$ be base-monotone, i-slice-continuous. If $H_1,H_2\subset(0,+\infty),\ H_1\cap H_2=\emptyset,\ {\cal{K}}(H_1)={\cal{K}}(H_2)=+\infty$ then ${\cal{K}}(H_1\cup H_2)=+\infty$.
\end{prp}
\P By i-slice-continuity we know that $\lim\limits_{x\to+\infty}{\cal{K}}(H_1^{x-})={\cal{K}}(H_1)$ and similarly for $H_2$ and $H_1\cup H_2$. Base-monotonicity gives that 
$$\min\{{\cal{K}}(H_1^{x-}),{\cal{K}}(H_2^{x-})\}\leq{\cal{K}}(H_1^{x-}\cup H_2^{x-})={\cal{K}}((H_1\cup H_2)^{x-}).$$ 
The limit of the left hand side is infinite hence so is the limit of the right hand side which gives the statement.
\Pes

\begin{prp}\label{pufmm}Let ${\cal{K}}$ be base-monotone, i-slice-continuous, subset-finite. If $H_1,H_2\subset(0,+\infty),\ {\cal{K}}(H_1)<\infty,\ {\cal{K}}(H_2)<\infty$ then ${\cal{K}}(H_1\cup H_2)<\infty$.
\end{prp}
\P By subset-finiteness we can assume that $H_1,H_2$ are disjoint. By i-slice-continuity we know that $\lim\limits_{x\to+\infty}{\cal{K}}(H_1^{x-})={\cal{K}}(H_1)$ and similarly for $H_2$ and $H_1\cup H_2$. Base-monotonicity gives that 
$${\cal{K}}((H_1\cup H_2)^{x-})={\cal{K}}(H_1^{x-}\cup H_2^{x-})\leq\max\{{\cal{K}}(H_1^{x-}),{\cal{K}}(H_2^{x-})\}.$$ 
The limit of the right hand side is finite hence so is the limit of the left hand side which gives the statement.
\Pes

\begin{prp}\label{peulekf}Let ${\cal{K}}$ be part-slice-continuous, i-slice-continuous, finite-independent. Moreover let the finite intervals be in $Dom\ {\cal{K}}$. Then $\forall\epsilon>0$ there is a unbounded $H\in Dom\ {\cal{K}}$ such that $\lambda(H)\leq\epsilon$ and ${\cal{K}}(H)<\infty$.
\end{prp}
\P Let $I_0=[0,\frac{\epsilon}{2}]$. If we have choosen $I_0,\dots,I_{n-1}$ already then find interval $I_n$ such that $I\subset[n,n+\frac{\epsilon}{2^{n+1}}]$ and ${\cal{K}}(\bigcup\limits_{i=0}^n I_i)<1$. Such interval can be found because consider $H_n=\bigcup\limits_{i=0}^n I_i\cup[n,n+\frac{\epsilon}{2^{n+1}}]$. If ${\cal{K}}(H_n)\geq 1$ then let $f(x)={\cal{K}}\big(\bigcup\limits_{i=0}^n I_i\cup[n,x]\big)\ (n\leq x\leq n+\frac{\epsilon}{2^{n+1}})$. By part-slice-continuity $f$ is continuous and by finite-independece $f(n)={\cal{K}}(\bigcup\limits_{i=0}^n I_i)<1$. Hence there is $x>n$ such that $f(x)<1$ still. Let $I_n=[n,x]$.

Let $H=\bigcup\limits_{i=0}^{\infty} I_i$. Clearly $H$ is unbounded, $\lambda(H)\leq\epsilon$ and by i-slice-continuity ${\cal{K}}(H)\leq 1$.
\Pes

\begin{prp}Let ${\cal{K}}$ be part-slice-continuous, i-slice-continuous, finite-independent, base-monotone, subset-finite and ${\cal{K}}([0,+\infty))=+\infty,\ {\cal{K}}((-\infty,0])=-\infty$. Moreover let the bounded sets be small for sets with infinite mean and let the finite intervals be in $Dom\ {\cal{K}}$. Then $\forall h\in\mathbb{R}$ there is $H\in Dom\ {\cal{K}}$ such that both $H^{0-},H^{0+}$ are unbounded and ${\cal{K}}(H)=h$.
\end{prp}
\P Let $h\in\mathbb{R}$.

First let us observe that base-monotonicity gives that $\forall a\in\mathbb{R}^+\ {\cal{K}}([a,+\infty))=+\infty,\ {\cal{K}}((-\infty,-a])=-\infty$ since base-monotonicity yields that $+\infty={\cal{K}}([0,+\infty))={\cal{K}}([0,a)\cup[a,+\infty))\leq\max\{{\cal{K}}([0,a)),{\cal{K}}([a,+\infty))\}$.

According to \ref{peulekf} we can construct unbounded $H_1,H_2\subset\mathbb{R}$ such that $H_1\subset(0,+\infty),\ H_2\subset(-\infty,0)$ and ${\cal{K}}(H_1)<\infty,\ {\cal{K}}(H_2)<\infty$. By \ref{pufmm} $k={\cal{K}}(H_1\cup H_2)<\infty$. If $k=h$ then we are done. 

Say $k<h$ (the other inequality is similar). Let $$f(x)={\cal{K}}\big(H_1\cup H_2\cup [h,x)\big)\ \ (x\in[h,+\infty)).$$ 
By part-slice-continuity $f$ is continuous. By \ref{pafaiui}, our first observation and i-slice-continuity we get that $\lim\limits_{x\to+\infty}f(x)=+\infty$. Hence there is $x\in[h,+\infty)$ such that ${\cal{K}}(H_1\cup H_2\cup [h,x))=h$.
\Pes

\begin{prp}If ${\cal{K}}$ is monotone and ${\cal{K}}(H)=+\infty,\ x\in\mathbb{R}$ then ${\cal{K}}(H^{+x})=+\infty$.
\end{prp}
\P Suppose that ${\cal{K}}(H^{x+})$ was finite for some $x\in\mathbb{R}$. Clearly $\sup (H\cap(-\infty,x))\leq\inf H^{x+}$ holds which gives that ${\cal{K}}(H)={\cal{K}}\big((H\cap(-\infty,x))\cup H^{x+}\big)\leq{\cal{K}}(H^{x+})<\infty$ that is a contradiction.
\Pes

\begin{prp}If ${\cal{K}}$ is monotone, bounded-finite and not finite then ${\cal{K}}$ is not Cantor-continuous.
\end{prp}
\P If ${\cal{K}}$ is not finite then there is $H\in Dom\ {\cal{K}}$ such that ${\cal{K}}(H)=\pm\infty$, say ${\cal{K}}(H)=+\infty$.
Find an $x_0\in\mathbb{R}$ such that $H^{x_0-}\ne\emptyset$. Let $C_0=H^{x_0-},\ C_n=C_0\cup H^{(x_0+n)+}$. Evidently $\bigcap\limits_{i=1}^{\infty}C_i=C_0$. Then $\forall i>0\ {\cal{K}}(C_i)=+\infty$ because otherwise by bounded-finiteness we would get that ${\cal{K}}(H)={\cal{K}}(C_i\cup (H\cap(x_0,x_0+i))<+\infty$. But ${\cal{K}}(C_0)<+\infty$ showing that ${\cal{K}}$ is not Cantor-continuous.
\Pes

\begin{lem}\label{lsfs}Let ${\cal{K}}$ be strong-base-monotone, i-slice-continuous. If $H,K\in Dom\ {\cal{K}},\ H\cap K=\emptyset,\ {\cal{K}}(H)=+\infty,\ K$ is bounded and ${\cal{K}}(H\cup K)\ne{\cal{K}}(K)$ then $K$ small for $H$.
\end{lem}
\P By i-slice-continuity $\lim\limits_{x\to+\infty}{\cal{K}}(H^{x-})=+\infty$. By base-monotonicity there are constants $0\leq c_x,d_x\leq 1$ such that $c_x+d_x=1$ and 
$${\cal{K}}(H^{x-}\cup K)=c_x{\cal{K}}(H^{x-})+d_x{\cal{K}}(K).$$ 
By assumption there is $x\in\mathbb{R}$ such that $c_x>0$ since otherwise we would get that ${\cal{K}}(H\cup K)={\cal{K}}(K)$. Strong-base-monotonicity implies that if $x<y$ then $c_x\leq c_y$. But then we get that $\lim\limits_{x\to+\infty}{\cal{K}}(H^{x-}\cup K)=+\infty$.
\Pes

\begin{prp}\label{pbss}Let ${\cal{K}}$ be strong-base-monotone, i-slice-continuous and $H,K\in Dom\ {\cal{K}},\ H\cap K=\emptyset,\ {\cal{K}}(H)=+\infty,\ K$ being bounded implies that ${\cal{K}}(H\cup K)\ne{\cal{K}}(K)$.
Then the bounded sets are small for sets with infinite mean.
\end{prp}
\P Let $H,K\in Dom\ {\cal{K}}$ such that $K$ is bounded and ${\cal{K}}(H)=+\infty$. We can assume that $H,K$ are disjoint. Then by \ref{lsfs} we get the statement.
\Pes

\begin{ex}Similarly to the proof of \ref{piscifhi} one can show that for every $0\leq s\leq 1,\ \epsilon>0$ there is $H\subset\mathbb{R}$ such that $\mu^s(H)<\epsilon$ and $Avg^s(H)=+\infty$.
\end{ex}

\begin{ex}${\cal{M}}^{\mu}$ is strong-base-monotone.
\end{ex}
\P We have to show the "strong" part only. Set ${\cal{K}}={\cal{M}}^{\mu}$.
Let $H_1,H_2,K\in Dom\ {\cal{K}}$ be bounded sets such that $H_1\subset H_2,\ H_2\cap K=\emptyset$. Then 
$${\cal{K}}(H_1\cup K)=\frac{\mu(H_1){\cal{K}}(H_1)+\mu(K){\cal{K}}(K)}{\mu(H_1)+\mu(K)}$$
and similarly for $H_2$. We have to show that 
$$\frac{\mu(H_1)}{\mu(H_1)+\mu(K)}\leq\frac{\mu(H_2)}{\mu(H_2)+\mu(K)}$$
which is straightforward.
\Pes

\begin{ex}The bounded sets are small for sets with infinite mean with respect to ${\cal{M}}^{\mu}$. 
\end{ex}
\P Apply \ref{pbss}.
\Pes

\begin{ex}For $Avg$ there are $H,K\in Dom\ Avg$ such that $K$ is bounded, $Avg(H)=+\infty$ and $Avg(H\cup K)<\infty$.

Let $K=[0,1]$ and $H$ be a $0.9$-set such that $Avg(H)=+\infty$. Then $Avg(H\cup K)=Avg(K)=0.5$.
\end{ex}

\begin{prp}${\cal{K}}$ is subset-finite iff $|{\cal{K}}(K)|=+\infty,\ K\subset H$ implies that either $H\not\in Dom\ {\cal{K}}$ or $|{\cal{K}}(H)|=+\infty$. \Pes
\end{prp}

\begin{cor}If ${\cal{K}}$ is subset-finite then the bounded sets are small for sets with infinite mean. \Pes
\end{cor}

\begin{ex}${\cal{M}}^{\mu}$ is subset-finite, bounded-finite.
\end{ex}
\P Both follow from 
\begin{equation}\label{eq1}
{\cal{M}}^{\mu}(A\cup^* B)=\frac{\mu(A){\cal{M}}^{\mu}(A)+\mu(B){\cal{M}}^{\mu}(B)}{\mu(A)+\mu(B)}
\end{equation}
when $A,B\in Dom\ {\cal{M}}^{\mu}$ (see \cite{lambm} Proposition 2.8).

If $|{\cal{M}}^{\mu}(H)|<\infty,\ K\subset H$ then set $A=K,\ B=H-K$ and then (\ref{eq1}) gives that ${\cal{M}}^{\mu}$ is subset-finite.

If $|{\cal{M}}^{\mu}(H)|<\infty,\ K\in Dom\ {\cal{M}}^{\mu}$ is bounded then set $A=H,\ B=K-H$ and then (\ref{eq1}) gives that ${\cal{M}}^{\mu}$ is bounded finite.
\Pes

\begin{ex}${\cal{M}}^{\mu}$ is interval-continuous.
\end{ex}
\P It simply follows from (\ref{eq1}) if we substitute $B$ with $I+x$ and use the fact that $\mu$ is $\epsilon-\delta$ absolut continuous with respect to $\lambda$.
\Pes

\begin{ex}$Avg^1$ is interval-infinite.
\end{ex}
\P It can be derived from (\ref{eq1}) if we substitute $B$ with $I+x$ and remark that $\lambda(I+x)=\lambda(I)$ and $Avg^1(I+x)=Avg^1(I)+x$.
\Pes

\begin{prp}If ${\cal{K}}$ is monotone then ${\cal{K}}(H^{+x})$ is increasing ($H\in Dom\ {\cal{K}}$).
\end{prp}
\P Let $x<y\ (x,y\in\mathbb{R})$. Let $H_1=H\cap[x,y],\ H_2=H\cap[y,+\infty)$. Then $\sup H_1\leq\inf H_2$ which implies that ${\cal{K}}(H^{+x})={\cal{K}}(H_1\cup H_2)\leq{\cal{K}}(H_2)={\cal{K}}(H^{+y})$.
\Pes

\begin{cor}If ${\cal{K}}$ is monotone, $x<y$ and ${\cal{K}}(H^{+x})=+\infty$ then ${\cal{K}}(H^{+y})=+\infty$. If $x<y$ and ${\cal{K}}(H^{+y})<+\infty$ then ${\cal{K}}(H^{+x})<+\infty.$ \Pes
\end{cor}

\begin{prp}Let ${\cal{K}}$ be shift-invariant, $H\subset\mathbb{R},\ x\in\mathbb{R},\ x\ne 0$ such that $H+x=H$. Then $H\not\in Dom\ {\cal{K}}$.
\end{prp}
\P Clearly it would mean that ${\cal{K}}(H)+x={\cal{K}}(H)$.
\Pe

The following can be shown similarly.

\begin{prp}Let ${\cal{K}}$ be homogeneous, $H\subset\mathbb{R},\ \alpha\in\mathbb{R}$ such that $\alpha H=H$. Then ${\cal{K}}(H)=0$ or $\alpha=1$. \Pes
\end{prp}

We close this section with some examples.

\begin{ex}Let  ${\cal{K}}=Avg$. Let $H_1\subset[0,1]$ be a set such that $H_1^{(1-\frac{1}{n})+}$ has Hausdorff dimension $\frac{1}{2}$ and $\mu^{\frac{1}{2}}(H_1^{(1-\frac{1}{n})+})>0$ for all $n\in\mathbb{N}$.
Let $H_2\subset[2,+\infty)$ be a set with Hausdorff dimension $\frac{1}{3}$ such that $Avg(H_2)=+\infty$. 
Then ${\cal{K}}(H^{x+})<1$ when $x<1$ and if $x\geq 1$ implies that ${\cal{K}}(H^{x+})=+\infty$. 
\end{ex}

\begin{ex}We present an example that is not interval-infinite.

Let us take the Borel measure associated to the harmonic mean: $\mu([a,b])=\mu_f([a,b])=\frac{1}{a^2}-\frac{1}{b^2}\ (a,b>0)$ and the mean by that measure for $H$ with $0<\mu(H)<+\infty$
$${\cal{M}}^{\mu}(H)=\frac{\int\limits_Hx d\mu}{\mu(H)}$$ 
(see 3.7 in \cite{lambm}). Take $H=[1,2],I=[1,2]$.
Then  $${\cal{M}}^{\mu}(H\cup I+x)=\frac{\big(\frac{1}{1}-\frac{1}{2}\big)+\big(\frac{1}{1+x}-\frac{1}{2+x}\big)}{\big(\frac{1}{1^2}-\frac{1}{2^2}\big)+\big(\frac{1}{(1+x)^2}-\frac{1}{(2+x)^2}\big)}\to\frac{\frac{1}{1}-\frac{1}{2}}{\frac{1}{1^2}-\frac{1}{2^2}}$$
when $x$ tends to infinity.

\smallskip

It is also an example for a finite mean. For simplicity let us restrict ${\cal{M}}^{\mu}$ for measurable subsets of $[1,+\infty)$.

First let us show that ${\cal{M}}^{\mu}([1,+\infty))=2$. Easy calculation shows that
$${\cal{M}}^{\mu}([1,+\infty))=\frac{\int\limits_1^{\infty}\frac{2}{x^2}\cdot d\lambda}{1}=2\bigg[-\frac{1}{x}\bigg]_1^{\infty}=2.$$
Now if $K\subset[1,+\infty),\ 0<\mu(K)<+\infty$ then $\int\limits_K xd\mu\leq\int\limits_{[1,+\infty)} xd\mu=2$ therefore ${\cal{M}}^{\mu}(K)$ is finite.

\smallskip

For similar reason this mean is not limit-finite: ${\cal{M}}^{\mu}([a,+\infty))=2a\ (a>0)$.
\end{ex}

\begin{ex}$Avg^1$ is not limit-finite. 

Let $H=\bigcup\limits_{i=1}^{\infty}[i,\frac{1}{2^i}]$. Then 
$$Avg^1(H)=\frac{1}{2}\frac{\sum\limits_{i=1}^{\infty}(i+\frac{1}{2^i})^2-i^2}{\sum\limits_{i=1}^{\infty}\frac{1}{2^i}}=
\frac{1}{2}\frac{\sum\limits_{i=1}^{\infty}2\frac{i}{2^i}+\frac{1}{2^{2i}}}{1}$$
that is clearly finite. For $n\in\mathbb{N}$ we get
$$Avg^1(H^{n+})=\frac{1}{2}\frac{\sum\limits_{i=n}^{\infty}2\frac{i}{2^i}+\frac{1}{2^{2i}}}{\sum\limits_{i=n}^{\infty}\frac{1}{2^i}}=\frac{1}{2}\frac{2\frac{n+1}{2^{n-1}}+\frac{1}{3\cdot 4^{n-1}}}{\frac{1}{2^{n-1}}}=n+1+\frac{1}{3\cdot 2^{n}}$$
showing that $Avg^1(H^{n+})-n\not\to 0$.
\end{ex}

\begin{ex}Let $\mu$ a Borel measure on $(0,+\infty)$ such that $\mu(\{x\})=0\ (x\in\mathbb{R})$. Then ${\cal{M}}^{\mu}$ is i-strong-internal.
\end{ex}
\P If ${\cal{K}}$ is not i-strong-internal then there is a set $H$ with $\sup H=+\infty$ and $H^{n+}$ consists of isolated points for some $n\in\mathbb{N}$ such that ${\cal{K}}(H)=+\infty$. It is trivial that for such $H$ we get that ${\cal{M}}^{\mu}(H)\leq n$.
\Pes

%------------------------------------------------------------------------------------------------------------------------------Examples---------------------------
\section{Examples}

In this section we present some examples on means that are defined on some unbounded sets as well.

\medskip 

\begin{ex}Let $H\subset[0,+\infty),\ H\ne\emptyset$. Let ${\cal{K}}^0(H)=\inf H$ and let
Let $\beta<\omega_1$ be an ordinal number and ${\cal{K}}^{\alpha}$ be already defined for $\alpha<\beta$. 
If $\beta$ is a successor ordinal, $\beta=\alpha+1$ then let 
$${\cal{K}}^{\beta}(H)=
\begin{cases}
\inf[{\cal{K}}^{\alpha}(H),+\infty)-H&\text{if }[{\cal{K}}^{\alpha}(H),+\infty)-H\ne\emptyset\\
{\cal{K}}^{\alpha}(H)&\text{otherwise.}
\end{cases}
$$
If $\beta<\omega_1$ is a limit ordinal then set ${\cal{K}}^{\beta}(H)=\sup\{{\cal{K}}^{\alpha}(H):\alpha<\beta\}$.
\end{ex}

\begin{ex}Let $H\subset[1,+\infty)$ be unbounded. Let $${\cal{K}}(H)=\inf H+\sum\limits_{i=1}^{\infty}\frac{f(i)}{i}$$ where 
$$f(i)=
\begin{cases}
1&\text{if }[i,i+1)\cap H\ne\emptyset\\
0&\text{otherwise.}
\end{cases}
$$
\end{ex}

\begin{ex}Let $(a_n)$ be an increasing sequence such that $a_n\to+\infty$. Let $H\subset[0,+\infty)$. Let
$${\cal{K}}(H)=
\begin{cases}
\sup\{a_n:a_n\in H\}&\text{if }\exists n\ a_n\in H\\
\inf H&\text{otherwise.}
\end{cases}
$$\end{ex}

\begin{ex}Let ${\cal{K}}$ be a mean defined on bounded subsets. Let us extend ${\cal{K}}$ in the simpliest way: 
$$\tilde{{\cal{K}}}(H)=
\begin{cases}
+\infty&\text{if }H\text{ is unbounded}\\
{\cal{K}}(H)&\text{otherwise.}
\end{cases}
$$
\end{ex}

\begin{ex}If $H\subset(0,+\infty)$ then set $\frac{1}{H}=\{\frac{1}{h}:h\in H\}$. Let ${\cal{K}}$ be a mean defined on bounded subsets of $(0,+\infty)$ such that $H\in Dom\ {\cal{K}}$ implies that $\frac{1}{H}\in Dom\ {\cal{K}}$. If $\sup H=+\infty,\ \inf H>0$ then we can extend ${\cal{K}}$ to $H$ in the following way: $${\cal{K}}(H)=\frac{1}{{\cal{K}}\big(\frac{1}{H}\big)}.$$
\end{ex}

\begin{ex}We can slightly generalize the definition of mean-set ${\cal{MS}}^{hf}(H)$ given in \cite{lamis} Definition 15 .
For $0<\lambda(H)<+\infty$ let ${\cal{MS}}^{hf}(H)=\{x:\lambda(H^{-x})=\lambda(H^{+x})\}$. Clearly ${\cal{MS}}^{hf}(H)$ is a closed finite interval.
\end{ex}

%------------------------------------------------------------------------------------------------------------------------------Extending means---------------------------
\section{Extending means}

In this section we are going to investigate how the domain of a mean can be extended to unbounded sets as well. 

\begin{prp}\label{pembb}Let ${\cal{K}}$ be a monotone mean whose domain contains bounded sets only. If $H\subset\mathbb{R}$ such that $\forall x\in\mathbb{R}\ H^{x+},H^{x-}\in Dom\ {\cal{K}}$, $\inf H>-\infty$ then $\lim\limits_{x\to+\infty}{\cal{K}}(H^{x-})$ exists.
\end{prp}
\P Obviously ${\cal{K}}(H^{x-})$ is increasing.
\Pe

One could formulate a similar statement for $+\infty$.

\smallskip

Now we are going to define a way of extending means.

\begin{df}Let ${\cal{K}}$ be a mean whose domain contains bounded sets only. Let $H\subset\mathbb{R}$ such that $\forall x\in\mathbb{R}\ H^{x+},H^{x-}\in Dom\ {\cal{K}}$.
Set $$\hat{\cal{K}}(H)=\mathop{\lim\limits_{x\to-\infty}}_{y\to+\infty}{\cal{K}}(H\cap[x,y])$$ if the limit exits. 
\end{df}

\begin{prp}Evidently $Dom\ {\cal{K}}\subset Dom\ \hat{\cal{K}}$ , moreover if $H\in Dom\ {\cal{K}}$ then $\hat{\cal{K}}(H)={\cal{K}}(H)$ i.e. $\hat{\cal{K}}$ is an extension of ${\cal{K}}$. \Pes
\end{prp}

\begin{rem}By \ref{pembb} if $H\subset\mathbb{R}$ unbounded such that $\forall x\in\mathbb{R}\ H^{x+},H^{x-}\in Dom\ {\cal{K}}$, $\inf H>-\infty$ or $\sup H<+\infty$ then ${\cal{K}}$ can be extended to $H$.
\end{rem}

\begin{prp}\label{pmls}$\hat{\cal{K}}(H)=h\in\bar{\mathbb{R}}$ iff for all sequences $(x_n),(y_n)$  such that $x_n\to-\infty,y_n\to+\infty$ $\lim\limits_{n\to\infty}{\cal{K}}(H\cap[x_n,y_n])=h$ holds. \Pes
\end{prp}

\begin{prp}Let $H\subset\mathbb{R}$ such that $\forall x\in\mathbb{R}\ H^{x+},H^{x-}\in Dom\ {\cal{K}}$. Then ${\cal{K}}$ can be extended to $H$ iff for all sequences $(x_n),(y_n)$  such that $x_n\to-\infty,y_n\to+\infty$ the limit $\lim\limits_{n\to\infty}{\cal{K}}(H\cap[x_n,y_n])$ always exists. \Pes
\end{prp}
\P Merging two sequences shows that both sequences have to provide the same limit hence \ref{pmls} is applicable.
\Pes

\begin{prp}Let ${\cal{K}}$ be a mean whose domain contains bounded sets only. Let $H\subset\mathbb{R}$ such that $\forall x\in\mathbb{R}\ H^{x+},H^{x-}\in Dom\ {\cal{K}}$. Let $f(x,y)={\cal{K}}(H\cap[x,y])$ be integrable.
If $\hat{\cal{K}}(H)$ exists then 
$$\hat{\cal{K}}(H)=\lim\limits_{p\to+\infty}\frac{1}{(2p)^2}\int\limits_{-p}^p\int\limits_{-p}^p{\cal{K}}(H\cap[x,y])dxdy.$$
\end{prp}
\P Denote $A=\hat{\cal{K}}(H)$. Let $A\in\mathbb{R}$. The cases $A=\pm\infty$ can be handled similarly.

Let $\epsilon>0$ be given. Choose $N\in\mathbb{R}^+$ such that if $x<-N<N<y$ then $|{\cal{K}}(H\cap[x,y])-A|<\frac{\epsilon}{3}$. Choose $M\in\mathbb{R}^+,\ M>N$ such that 
$$\frac{N(2N)^2}{(2M)^2}<\frac{\epsilon}{2}\text{ and }A-\frac{\epsilon}{2}<\frac{(A-\frac{\epsilon}{3})((2M)^2-(2N)^2)}{(2M)^2}.$$
Let $p>M$. 
Let us use the notations $$L=[-N,N]\times[-N,N],\ K=[-p,p]\times[-p,p]-L.$$ Then 
$$\lim\limits_{p\to+\infty}\frac{1}{(2p)^2}\int\limits_{-p}^p\int\limits_{-p}^p{\cal{K}}(H\cap[x,y])dxdy=$$
$$\lim\limits_{p\to+\infty}\frac{1}{(2p)^2}\int\limits_{L}{\cal{K}}(H\cap[x,y])dxdy+\frac{1}{(2p)^2}\int\limits_{K}{\cal{K}}(H\cap[x,y])dxdy.$$
Clearly $x,y\in L$ implies that $|{\cal{K}}(H\cap[x,y])|<N$ hence $$\left\vert\frac{1}{(2p)^2}\int\limits_{L}{\cal{K}}(H\cap[x,y])dxdy\right\vert<\frac{N(2N)^2}{(2M)^2}<\frac{\epsilon}{2}.$$
Moreover
$$A-\frac{\epsilon}{2}<\frac{(A-\frac{\epsilon}{3})((2p)^2-(2N)^2)}{(2p)^2}=\frac{1}{(2p)^2}\int\limits_{K}A-\frac{\epsilon}{3}dxdy<\frac{1}{(2p)^2}\int\limits_{K}{\cal{K}}(H\cap[x,y])dxdy$$
$$<\frac{1}{(2p)^2}\int\limits_{K}A+\frac{\epsilon}{3}dxdy=\frac{(A+\frac{\epsilon}{2})((2p)^2-(2N)^2)}{(2p)^2}<A+\frac{\epsilon}{2}.$$
Therefore if $p>M$ then $$A-\epsilon<\frac{1}{(2p)^2}\int\limits_{-p}^p\int\limits_{-p}^p{\cal{K}}(H\cap[x,y])dxdy<A+\epsilon.$$
\Pes

\begin{prp}If ${\cal{K}}_1,{\cal{K}}_2$ are two means whose domain contains bounded sets only, $Dom\ {\cal{K}}_1=Dom\ {\cal{K}}_2$ and ${\cal{K}}_1\leq{\cal{K}}_2$ then $\hat{\cal{K}}_1\leq\hat{\cal{K}}_2$.  \Pes
\end{prp}

For some means one can ask if the straightforward (algebric) generalization of the mean to unbounded sets equals to the extension that we have just defined. We investigate two means in this respect.  

\begin{prp}Let ${\cal{K}}={\cal{M}}^{\mu}$ restricted to the bounded sets. If $0<\mu(H)<+\infty$ and $\hat{\cal{K}}(H)$ exists then $\hat{\cal{K}}(H)={\cal{M}}^{\mu}(H)$.
\end{prp}
\P Clearly $$\hat{\cal{K}}(H)=\mathop{\lim\limits_{x\to-\infty}}_{y\to+\infty}\frac{\int\limits_{H\cap[x,y]} zd\mu(z)}{\mu(H\cap[x,y])}=\frac{\int\limits_{H} zd\mu(z)}{\mu(H)}={\cal{M}}^{\mu}(H).$$
\Pes

%\begin{lem}Let $\inf H>-\infty,\ 0<\mu(H)<+\infty,\ {\cal{M}}^{\mu}(H)=+\infty$. Then $\forall x\in\mathbb{R}\ {\cal{M}}^{\mu}(H^{x+})=+\infty$.
%\end{lem}
%\P Set $A=\sup\{y: \mu(H^{y+})=\mu(H)\}$. If $x\leq A$ then the statement is obvious. If $A<x$ then it can be seen from 
%$${\cal{M}}^{\mu}(H)=\frac{\mu(H^{x-}){\cal{M}}^{\mu}(H^{x-})+\mu(H^{x+}){\cal{M}}^{\mu}(H^{x+})}{\mu(H)}.$$
%\Pes

\begin{ex}Let ${\cal{K}}={\cal{M}}^{lis}$ (i.e. ${\cal{K}}(H)=\frac{\varliminf H+\varlimsup H}{2}$ for a bounded $H$). Then $\hat{\cal{K}}(H)$ is finite iff there is $n\in\mathbb{R}^+$ such that there is no finite accumulation point of $H-[-n,n]$. We can conclude that $\hat{\cal{K}}(H)\ne\frac{\varliminf H+\varlimsup H}{2}$ in general (e.g. for $H=\mathbb{N}\cup\{\frac{1}{n}:n\in\mathbb{N}\},\ \hat{\cal{K}}(H)=0,\ \frac{\varliminf H+\varlimsup H}{2}=+\infty$). It is also clear that $H\in\ Dom\ \hat{\cal{K}} \iff |H''\cap\{-\infty,+\infty\}|<2$.
\Pes
\end{ex}

%-----------------------------------------------------------------------------------------------------------------------------Inherited properties---------------------------
\subsection{Inherited properties}

In this section we investigate some properties which the extension inherite from the original mean.

\begin{prp}If ${\cal{K}}$ is monotone then so is $\hat{\cal{K}}$.
\end{prp}
\P Let $H_1,H_2,H_1\cup H_2\in Dom\ \hat{\cal{K}}$ and let $\sup H_1\leq\inf H_2$. If $x<\sup H_1$ and $\inf H_2<y$ then $\sup H_1\cap[x,y]\leq\inf H_2\cap[x,y]$ which gives that
$${\cal{K}}(H_1\cap[x,y])\leq{\cal{K}}((H_1\cup H_2)\cap[x,y])\leq{\cal{K}}(H_2\cap[x,y])$$
which gives the statement when $x\to-\infty,y\to+\infty$.
\Pes

\begin{prp}If ${\cal{K}}$ is base-monotone then so is $\hat{\cal{K}}$.
\end{prp}
\P Let $H_1,H_2,H_1\cup H_2\in Dom\ \hat{\cal{K}},\ H_1\cap H_2=\emptyset$. If $x,y\in\mathbb{R}$ then $(H_1\cap[x,y])\cap (H_2\cap[x,y])=\emptyset$ which gives that
$$\min\{{\cal{K}}(H_1\cap[x,y]),{\cal{K}}(H_2\cap[x,y])\}\leq{\cal{K}}((H_1\cup H_2)\cap[x,y])\leq\max\{{\cal{K}}(H_1\cap[x,y]),{\cal{K}}(H_2\cap[x,y])\}$$
which gives the statement when $x\to-\infty,y\to+\infty$.
\Pes

\begin{prp}If ${\cal{K}}$ is monotone, symmetric then $\hat{\cal{K}}$ is symmetric as well.
\end{prp}
\P Let $H\in Dom({\cal{K}})$ symmetric i.e. $\exists s\in\mathbb{R}\ T_s(H)=H$ where $T_s$ denote the reflection to point $s\in\mathbb{R}$ that is $T_s(x)=2s-x\ (x\in\mathbb{R})$.
Let $\hat{\cal{K}}(H)=A\in\bar{\mathbb{R}}$. Suppose that $s\ne A$. Take a neighbourhood $K$ of $A$ such that $s\not\in K$. We know that there are numbers $N,M$ such that $x\leq N<M\leq y$ implies that ${\cal{K}}(H\cap[x,y])\in K$. We can assume that $N<s,\ M=T_s(N)$. Then clearly ${\cal{K}}(H\cap[N,M])=s$ which is a contradiction.
\Pe

By the defnition it is clear that

\begin{prp}If ${\cal{K}}$ is slice-continuous then $\hat{\cal{K}}$ is i-slice-continuous. \Pes
\end{prp}

\begin{prp}If ${\cal{K}}$ is shift-invariant then so is $\hat{\cal{K}}$. 
\end{prp}
\P The obvious facts that $H^{x+}=(H+y)^{(x+y)+},(H+y)^{x+}=H^{(x-y)+}$ and similarly for $H^{x-}$ give the statement.
\Pes

\begin{prp}If ${\cal{K}}$ is disjoint-monotone then so is $\hat{\cal{K}}$. 
\end{prp}
\P Let $H_1,H_2,H_1\cup H_2\in Dom({\cal{K}}),\ H_1\cap H_2=\emptyset$. Let $\hat{\cal{K}}(H_1)\leq\hat{\cal{K}}(H_2)$. 

First suppose that $\hat{\cal{K}}(H_1)<\hat{\cal{K}}(H_2)$. Then there is $N\in\mathbb{R}$ such that if $x<-N<N<y$ then ${\cal{K}}(H_1\cap[x,y])\leq{\cal{K}}(H_2\cap[x,y])$. 
Then by disjoint-monotonicity of  ${\cal{K}}$ we get that ${\cal{K}}(H_1\cap[x,y])\leq{\cal{K}}((H_1\cup H_2)\cap[x,y])\leq{\cal{K}}(H_2\cap[x,y])$. If we take the limit we get the statement.

Now let $\hat{\cal{K}}(H_1)=\hat{\cal{K}}(H_2)\in\mathbb{R}$. Take two sequences $(x_n),(y_n)$ such that $x_n\to-\infty,\ y_n\to+\infty$. Let $(x_n'),(y_n')$ be the subsequences of $(x_n),(y_n)$ for which ${\cal{K}}(H_1\cap[x_n',y_n'])\leq{\cal{K}}(H_2\cap[x_n',y_n'])$ holds. Let us denote he remainder of the original sequences by $(x_n''),(y_n'')$ for which ${\cal{K}}(H_1\cap[x_n'',y_n''])>{\cal{K}}(H_2\cap[x_n'',y_n''])$ holds. 

For the first sequences we have ${\cal{K}}(H_1\cap[x_n',y_n'])\leq{\cal{K}}((H_1\cup H_2)\cap[x_n',y_n'])\leq{\cal{K}}(H_2\cap[x_n',y_n'])$. If $n\to\infty$ then we get that  $\hat{\cal{K}}(H_1)=\hat{\cal{K}}(H_1\cup H_2)=\hat{\cal{K}}(H_2)$. Similarly using the second sequences we get exactly the same $\hat{\cal{K}}(H_1)=\hat{\cal{K}}(H_1\cup H_2)=\hat{\cal{K}}(H_2)$ that completes the proof.
\Pes

%-----------------------------------------------------------------------------------------------------------------------------Some properties of ${\cal{M}}^{\mu}$---------------------------
\section{Some properties of ${\cal{M}}^{\mu}$}

Here we examine ${\cal{M}}^{\mu}\ (Avg^1)$ more closely in the sense of infinite behavior.

\begin{prp}\label{paii}If $H\subset\mathbb{R},\ \inf H>-\infty,\ \mu(H)=+\infty$ then $\hat{\cal{K}}(H)=+\infty$ for ${\cal{K}}={\cal{M}}^{\mu}$.
\end{prp}
\P Let $n>\inf H$. Then 
$$\hat{\cal{K}}(H)=\lim\limits_{m>n,m\to+\infty}{\cal{M}}^{\mu}(H^{m-})=$$
$$\lim\limits_{m>n,m\to+\infty}\frac{\mu(H^{n-}){\cal{M}}^{\mu}(H^{n-})+\mu(H\cap[n,m]){\cal{M}}^{\mu}(H\cap[n,m])}{\mu(H^{n-})+\mu(H\cap[n,m])}\geq$$
$$\lim\limits_{m>n,m\to+\infty}\frac{\mu(H\cap[n,m])\cdot n}{\mu(H^{n-})+\mu(H\cap[n,m])}=$$
$$\lim\limits_{m>n,m\to+\infty}\frac{n}{\frac{\mu(H^{n-})}{\mu(H\cap[n,m])}+1}=n$$
because $\mu(H\cap[n,m])\to+\infty$ which gives that $\frac{\mu(H^{n-})}{\mu(H\cap[n,m])}\to 0$.
\Pes

\begin{prp}Let ${\cal{K}}={\cal{M}}^{\mu}$. If $\hat{\cal{K}}(H^{0-})>-\infty,\ \hat{\cal{K}}(H^{0+})<+\infty$ then $\hat{\cal{K}}(H)$ exists and $|\hat{\cal{K}}(H)|<\infty$.
\end{prp}
\P We could simply refer to \ref{pufmm} but we give a direct proof too.

By \ref{paii} we get that $\mu(H^{0-})<\infty,\mu(H^{0+})<\infty$.
Clearly $$\hat{\cal{K}}(H)=\mathop{\lim\limits_{x\to-\infty}}_{y\to+\infty}\frac{\int\limits_{H\cap[x,y]} zd\mu(z)}{\mu(H\cap[x,y])}=
\mathop{\lim\limits_{x\to-\infty}}_{y\to+\infty}\frac{\int\limits_{H\cap[x,0]} zd\mu(z)+\int\limits_{H\cap[0,y]} zd\mu(z)}{\mu(H\cap[x,y])}=$$
$$\mathop{\lim\limits_{x\to-\infty}}_{y\to+\infty}\frac{\mu(H\cap[x,0]){\cal{M}}^{\mu}(H\cap[x,0])+\mu(H\cap[0,y]){\cal{M}}^{\mu}(H\cap[0,y])}{\mu(H\cap[x,y])}=$$
$$\frac{\mu(H^{0-}){\cal{M}}^{\mu}(H^{0-})+\mu(H^{0+}){\cal{M}}^{\mu}(H^{0+})}{\mu(H)}.$$
\Pes

\begin{ex}Let two sequences $(b_n),(c_n)$ be given such that $0<c_n<1,\ \sum\limits_{n=1}^{\infty}c_n<+\infty,\ b_n+c_n\leq b_{n+1}$. Let $I_n=[b_n,b_n+c_n],\ H=\bigcup\limits_{n=1}^{\infty}I_n$. Then $Avg^1(H)$ gets finite iff $\sum\limits_{n=1}^{\infty}b_n\cdot c_n<+\infty$. 
\end{ex}
\P Clearly
$$Avg^1(H)=\frac{1}{2}\frac{\sum\limits_{n=1}^{\infty}(b_n+c_n)^2-b_n^2}{\sum\limits_{n=1}^{\infty}c_n}=
\frac{\sum\limits_{n=1}^{\infty}c_n(b_n+\frac{1}{2}c_n)}{\sum\limits_{n=1}^{\infty}c_n}=\frac{\sum\limits_{n=1}^{\infty}b_n\cdot c_n+\frac{1}{2}\sum\limits_{n=1}^{\infty}c_n^2}{\sum\limits_{n=1}^{\infty}c_n}$$
which is finite iff $\sum\limits_{n=1}^{\infty}b_n\cdot c_n<+\infty$. 
\Pes

%---------------------------------------------------------------------------------------------------------------------------------------------------------thebibliography-----------------

{\footnotesize

%------------------------------------------------------------------------------------------------------------------------------------------------------------address---------------
\noindent
%Attila Losonczi,\\ Hungary 2120 Dunakeszi Bojtorj\'an u. 3/1\\

\noindent E-mail: alosonczi1@gmail.com\\

\end{document}